\documentclass[preprint,12pt]{elsarticle}

\usepackage{amsfonts,amsmath,amssymb,amsthm,amstext,txfonts}
\usepackage{extarrows,longtable,graphicx,hyperref,mathrsfs,url,booktabs}

\journal{ArXiv}

\newtheorem{theorem}{Theorem}[section]
\newtheorem{proposition}[theorem]{Proposition}
\newtheorem{lemma}[theorem]{Lemma}

\newtheorem{corollary}[theorem]{Corollary}
\numberwithin{equation}{section}
\allowdisplaybreaks

\begin{document}

\begin{frontmatter}
\title{Singular scaling properties of a tree-like fractal Dirichlet space}
\author{Guanhua Liu\fnref{L}}
\ead[L]{liu\_gh@tju.edu.cn}
\address{Center for Applied Mathematics and KL-AAGDM, Tianjin University, Tianjin, 300072, China}

\begin{abstract}
This work is a continuation of our recent work on a tree-like self-similar fractal equipped with a regular self-similar Dirichlet form. Effective resistance estimates are established on different combinations of cells, yielding a new formalism of exit time estimate. We characterize the weights of a self-similar measure such that the scaling property is disrupted, and obtain a first-step heat kernel estimate under weights of the opposite case. Simultaneously, an irregularly scaled capacity estimate is listed under all possible weights.
\end{abstract}

\end{frontmatter}

\section{Introduction}

\begin{figure}[htbp]
\centering\includegraphics[width=0.75\textwidth]{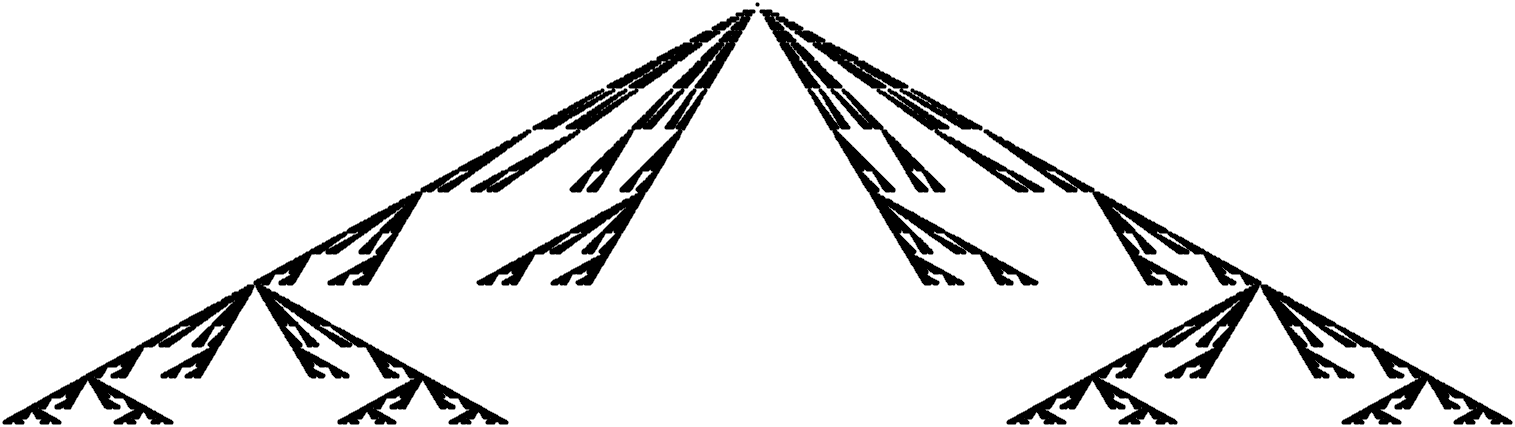}
\caption{The tree-like fractal}
\label{tree_fr}
\end{figure}

In a recent work \cite{HL}, we construct a self-similar fractal (see Figure \ref{tree_fr}) and corresponding resistance forms which constitute self-similar Dirichlet spaces with respect to properly selected measures. Compared with regular cases (cf.\ \cite{hk,L0} for post-critically finite fractals, \cite{bgk,ckw-e,ghl15,L1}, etc. for general metric measure spaces), on this fractal many singular properties are found, especially the separation between strong and weak elliptic Harnack inequalities. Particularly emphasized is some singular behaviour of the mean exit time function. We have shown in \cite[Theorem 3.9]{HL} that when the self-similar weights of the measure satisfy $w_2>w_0$, then the ratio between supremum and infimum of the mean exit time from some typical balls is unbounded. In this paper we focus on the remaining case $w_2\le w_0$. A strange type of scaling will be constructed, and for the particular case $w_0=w_2$, some first-step heat kernel estimates follow.

For self-similar fractals and resistance forms on them, we follow Kigami's notions \cite{ki}; while for operators on Dirichlet spaces, we follow notions given in the standard textbook \cite{fot} by Fukushima et. al.

\subsection{The tree-like fractal}
Let $K$ be the invariant set of the iterative function system 
\begin{align*}
F_0(x,y)=&\ \left(\frac{2}{9}x+\frac{8}{9\sqrt{3}}y,\frac{2}{3}y\right),\quad\quad F_2(x,y)=\left(\frac{1}{3}x-\frac{2}{3},\frac{1}{3}y-\frac{1}{\sqrt{3}}\right),\\
F_1(x,y)=&\ \left(\frac{2}{9}x-\frac{8}{9\sqrt{3}}y,\frac{2}{3}y\right),\quad\quad F_3(x,y)=\left(\frac{1}{3}x+\frac{2}{3},\frac{1}{3}y-\frac{1}{\sqrt{3}}\right)
\end{align*}
on $\mathbb{R}^2$, shown in Figure \ref{tree_fr}. For any $\omega=i_1\cdots i_n\in\mathbb{S}_n:=\{0,1,2,3\}^n$ (a word of length $n$), denote
$$F_\omega=F_{i_1}\circ\cdots\circ F_{i_n}\quad\mbox{and}\quad K_\omega=F_\omega(K).$$
Denote $i_1\cdots i_n$ as $i^n$ if every digit $i_j$ is identical to a same $i\in\{0,1,2,3\}$. Set
$$\mathbb{S}_0=\{\emptyset\},\quad F_\emptyset=\mathrm{id}\quad\mbox{and}\quad K_\emptyset=K.$$ Let $\mathbb{S}_*=\cup_{i=0}^\infty\mathbb{S}_n$ be the collection of all finite words.

On the contrary, let  $\mathbb{S}_\infty=\{0,1,2,3\}^\mathbb{N}$ be the collection of all infinite words, and $\pi:\mathbb{S}_\infty\to K$ be the projection, mapping every $\omega=i_1\cdots i_n\cdots\in\mathbb{S}_\infty$ to the unique point in $\cap_{n=1}^\infty K_{i_1\cdots i_n}$. For every $z\in K$, let $\omega_z$ be its first address, that is, the first element (in the reverse dictionary order) in $\pi^{-1}(\{z\})$.

Denote by $V_0$ the boundary of $K$, which consists of exactly three points
$$q_1=(0,0),\quad q_2=\bigg(-\frac{1}{2},-\frac{\sqrt{3}}{2}\bigg)\quad\mbox{and}\quad q_3=\bigg(\frac{1}{2},-\frac{\sqrt{3}}{2}\bigg),$$
respectively the invariant points of $F_0$ (also $F_1$), $F_2$ and $F_3$. Define inductively $V_{n+1}=\cup_{i=0}^3F_i(V_n)$, then $V_n\subset V_{n+1}$ for all $n\ge 0$. For convenience, we also denote
$$V_n^\uparrow:=\{F_\omega(q_1)\}_{\omega\in\mathbb{S}_n}\subset V_n.$$
It is easy to show both $V_\ast:=\cup_{n=0}^\infty V_n$ and $V_\ast^\uparrow:=\cup_{n=0}^\infty V_n^\uparrow$ are dense in $K$. Meanwhile, a standard Cantor set $\mathit{C}:=\pi\left(\{2,3\}^\mathbb{N}\right)$, the bottom of $K$, does not intersect $V_\ast^\uparrow$.

By \cite[Proposition 2.1]{HL}, there exists a unique resistance form $\mathcal{E}$ satisfying $$\mathcal{E}(u)=2\sum\limits_{i=0}^3\mathcal{E}(u\circ F_i)\quad\mbox{for all}\quad u\in C(K,R)$$
with effective resistance $R(q_1,q_2)=R(q_1,q_3)=1$. Further, every $F_i$ is a similitude under the resistance metric $(K,R)$ with ratio $\frac{1}{2}$, $R$ is geodesic and $\mathrm{diam}(K,R)=2$.

Given an open set $U\subset K$, we say $u\in C(K,R)$ is a harmonic function in $U$, if $\mathcal{E}(u)<\infty$ and $\mathcal{E}(u,\varphi)=0$ for every $\varphi\in C(K,R)$ supported in $U$ with $\mathcal{E}(\varphi)<\infty$. A harmonic function in $U$ is uniquely determined by its value on $\partial U$, while a globally harmonic function must be constant.

The following three harmonic functions (constructed in \cite[Lemmas 2.2, 2.3 and 2.5]{HL}) are frequently used in constructions:

\begin{itemize}
\item $u_{a_2,a_1,a_3}^-$, harmonic on $K\setminus V_0$, with boundary value $q_i\mapsto a_i$. It satisfies
$$u_{a_2,a_1,a_3}^-(x)=R(q_j,x)a_1+R(q_1,x)a_j\quad\mbox{on}\quad\overline{q_1q_j}\quad(j=2,3)$$
and is piecewise constant elsewhere. Further,
$$\mathcal{E}(u_{a_2,a_1,a_3}^-)=E_0(u|_{V_0})=(a_1-a_2)^2+(a_1-a_3)^2.$$
\item $u_\downarrow$, harmonic on $B(q_1,1)\setminus\{q_1\}$, with boundary value $q_1\mapsto 1$, $\mathit{C}\mapsto 0$. It satisfies
$$u_\downarrow(F_{2^m}(q_1))=4^{-m}\quad\mbox{for all}\quad m\in\mathbb{N}\quad\mbox{and}\quad\mathcal{E}(u_\downarrow)=3.$$
\item $u_\uparrow$, harmonic on $B(q_2,1)\setminus\{q_2\}$, with boundary value $q_2\mapsto 1$, $q_1\mapsto 0$, $F_{0^m23}(\mathit{C})\mapsto 0$ (for all $m\in\mathbb{N}$). It satisfies
$$u_\uparrow(F_{0^m2}(q_1))=4^{-1-m}\quad\mbox{and}\quad\mathcal{E}(u_\uparrow)=\frac{3}{2}.$$
\end{itemize}

Given every array of weights $\mathbf{w}:=(w_0,w_1,w_2,w_3)\in(0,1)^4$ with
$$w_0=w_1,\quad w_2=w_3\quad\mbox{and}\quad w_0+w_1+w_2+w_3=1,$$
there exists a unique self-similar measure $\mu=\mu_\mathbf{w}$ on $K$ such that $\mu(K)=1$ and
$$\mu(A)=\sum_{i=0}^3w_i\ \mu\left(F_i^{-1}(A)\right)\quad\mbox{for every Borel set}\quad A\subset K.$$
Note that $\mu(V_*)=0$, $\mu$ is never doubling, and $\mathcal{E}$ extends to a Dirichlet form on $L^2(K,\mu)$ (whose definition domain is denoted by $\mathcal{F}$) that is embedded into $C(K,R)$.

Note by \cite[Lemma 3.7]{HL} that for any $\mathbf{w}$,
\begin{equation}\label{int-psi}
\int_Ku_\uparrow d\mu\asymp\int_Ku_\downarrow d\mu\asymp 1.
\end{equation}
\subsection{Operators and isoperimetric conditions}

Fix a measure $\mu$ and consider the Dirichlet space $(K,R,\mu,\mathcal{E},\mathcal{F})$. For every open set $U\subset K$, let $\mathcal{F}(U)$ be the completion of $\mathcal{F}\cap C(K,R)$ under $\mathcal{E}(\cdot)+\|\cdot\|_{L^2(U,\mu)}$. The Green operator $G_\mu^U:L^2(U,\mu)\to\mathcal{F}(U)$ maps $f\in L^2(U,\mu)$ to the function $u\in\mathcal{F}(U)$ satisfying $\mathcal{E}(u,\varphi)=(f,\varphi)_{L^2(U,\mu)}$ for all $\varphi\in\mathcal{F}(U)$. Recall that, since $(K,R,\mathcal{E})$ is also a resistance form, then
$$G_\mu^U1(x)=R(x,U^c)\int_U\psi_U^xd\mu\quad\mbox{for all}\quad x\in U,$$
where $\psi_U^x$ is the function attaining $R(x,U^c)$, which itself is independent of $\mu$.

Closely related is the heat semigroup $\{P_t^U\}_{t>0}$ (simply written as $P_t$ if $U=K$), mapping every $f\in L^2(U)$ to the weak solution of homogeneous Dirichlet problem
$$\partial_t(u,\varphi)_{L^2(\mu)}+\mathcal{E}(u,\varphi)=0\quad\mbox{for all}\quad\varphi\in\mathcal{F}(U)$$
satisfying $\lim_{t\downarrow 0}\|u-f\|_{L^2(U,\mu)}=0$. If $P_t^U$ admits an integral kernel $p_t^U(x,y)$ for every $t>0$, we say $p_\cdot^U(\cdot,\cdot)$ is the $U$-local (global, if $U=K$) heat kernel of Dirichlet form $(\mathcal{E},\mathcal{F};L^2(\mu))$. Note that
$$G_\mu^Uf=\int_0^\infty P_t^Ufdt.$$
Functions $G_\mu^U1$ and $P_t^U1$ are respectively called the (mean) exit time and survival functions, according to what they stand for in the sense of stochastic processes.

We call any measurable function $W:K\times\mathbb{R}_+\to\mathbb{R}_+$ with $W(\cdot,0)\equiv 0$ as a quasi-scale, while call it a \textit{scale}, if $\lim\limits_{r\to\infty}W(x,r)=\infty$ and $W$ is weakly monotonic, that is, $W(x,R)\ge C^{-1}W(y,r)$ with some constant $C\in(0,\infty)$ whenever $R\ge r$ and $B(x,R)\supset B(y,r)$. Further, we say $W$ is a \textit{scaling function}, if there exist constants $C_W\ge 1$ and $0<\beta_1\le\beta_2<\infty$ such that
\begin{equation}\label{scal}
C_W^{-1}\left(\frac{R}{r}\right)^{\beta_1}\le\frac{W(x,R)}{W(y,r)}\le C_W\left(\frac{R}{r}\right)^{\beta_2}
\end{equation}
for all $0<r\le R<\infty$ and all $x,y\in M$ with $d(x,y)\le R$.

The scaling property (\ref{scal}) implies clearly that a scaling function is always a scale. Further, given a scaling function $W$, by considering $W'(x,r)=\int_0^r\frac{W(x,r')}{r'}dr'$ instead of $W$, we can always assume that $W(x,\cdot)$ is a homeomorphism on $[0,\infty)$ for every $x\in K$.

Since $\mathrm{diam}(K,R)=2$, all the properties of $(K,R,\mu,\mathcal{E},\mathcal{F})$ to be scaled concern balls of radius $r\le 2$ only. Thus it suffices to investigate $W(x,r)$ with $r\in[0,2]$, and the complete (quasi-)scale is easily constructed by letting $W(x,r)=(r/2)W(x,2)$ for all $r>2$. Further, when weak monotonicity is known, it even suffices to consider $W$ on $[0,1]$ (or $[0,r_0]$ with any fixed $r_0>0$).

Write $W(B):=W(x,r)$ for every ball $B=B(x,r)$. Given a scale $W$, we say the \emph{exit time estimate} $(\mathrm{E}_W)$ holds, if there are constants $C>0$ and $\delta,\sigma\in(0,1)$ such that for every ball $B=B(x_0,r)$ with $0<r<\mathrm{diam}(K,R)$, for $\mu$-a.e.\ $y\in B$,
\begin{align}
G_\mu^B1(y)\le CW(B)&\quad\quad\mbox{if}\quad\quad r<\sigma\cdot\mathrm{diam}(K,R);\tag*{$(\mathrm{E}_\le)$}\\
G_\mu^B1(y)\ge C^{-1}W(B)&\quad\quad\mbox{if}\quad\quad y\in\delta B.\tag*{$(\mathrm{E}_\ge)$}
\end{align}

We say the \textit{survival estimate} $(\mathrm{S}_W)$ holds, if there exist constants $\varepsilon,\delta\in(0,1)$ such that for every ball $B$ and all $0<t\le\delta W(B)$,
$$P_t^B1\ge\varepsilon\quad\mbox{on}\quad\delta B.$$

We say the \textit{capacity estimate} $(\mathrm{cap}_W)$ holds, if for any $\lambda\in(0,1)$, there exists a constant $C_\lambda>0$ such that for every ball $B$,
\begin{equation}\label{cap}
\mathrm{cap}_0(\lambda B,B):=\inf\bigg\{\mathcal{E}(u):u\in\mathcal{F}(B),0\le u\le 1,u|_{\lambda B}=1\bigg\}\le C_\lambda\frac{\mu(B)}{W(B)}.
\end{equation}

Recall that (cf.\ \cite{ghh+3}), these conditions are related by
\begin{equation}\label{imply}
(\mathrm{E}_W)\Rightarrow(\mathrm{S}_W)\Rightarrow(\mathrm{ABB}_W)\Rightarrow(\mathrm{cap}_W)
\end{equation}
with a same scaling function $W$. Here $(\mathrm{ABB}_W)$ stands for an artificial cutoff inequality (omitted here). It is not hard to check it holds also for general scales.

\subsection{A known estimate of exit time on $K$}

In the sequel, for every ball $B=B(q,r)$ with $q\in V_n$ and $0<r\le 2^{-n}$, set
\begin{align*}
B^\uparrow=&\ \cup\big\{B\cap K_\tau:\ \tau\in\mathbb{S}_n,q\in\{F_\tau(q_2),F_\tau(q_3)\}\big\};\\
B^\downarrow=&\ \cup\big\{B\cap K_\tau:\ \tau\in\mathbb{S}_n,q=F_\tau(q_1)\big\}.
\end{align*}
Clearly $B^\uparrow\cup B^\downarrow=B$ and $B^\uparrow\cap B^\downarrow=\{q\}$.

We have established a technique to construct and estimate $G_\mu^{B_n}1$ using self-similarity and properties of the basic harmonic functions $u_{a_2,a_1,a_3}^-$, $u_\uparrow$, $u_\downarrow$ on a sequence of typical balls $B_n=B(q_0,2^{-n})$, where $q_0=F_0(q_2)=F_2(q_0)$. The following lemma is fundamental:

\begin{lemma}\cite[Theorem 3.8]{HL}\label{exit}
For all $x\in B_n$,
$$G_\mu^{B_n}1(x)\asymp\begin{cases}
\left(2^nd_R(x,B_n^c)\right)^2\left(2^{-n}w_0^n+w_2^n\left(R(x,q_0)+2^{-2n}\right)\right),&\mbox{for}\quad x\in B_n^\uparrow;\\
\left(2^nd_R(x,B_n^c)\right)^2\left(2^{-n}w_0^n+2^{-2n}w_2^n\right),&\mbox{for}\quad x\in B_n^\downarrow.
\end{cases}$$
\end{lemma}

The same holds also for balls centred at $q'_0=F_0(q_3)$ by symmetry.

\section{Exit time for some typical regions}\label{S4}

In this section we use the same method as in \cite[Section 3]{HL} to obtain exit time estimates on some other typical regions, which will be used in Section \ref{S5} for the case $w_2\le w_0$.

The following notions are very important in the sequel. For every $q\in V_\ast^\uparrow$, let $n_\ast(q)\in\mathbb{N}$ be the least integer $n$ such that there exists $\omega\in\mathbb{S}_n$ satisfying $q=F_\omega(q_1)$. Note that either $\omega=\emptyset$ (so that $q=q_1$) or $\omega=\omega'i$ with some $\omega'\in\mathbb{S}_{n_\ast(q)-1}$ and $i\in\{2,3\}$ (so that $q=F_{\omega'}(q_0)$ or $q=F_{\omega'}(q'_0)$). In particular, this $\omega$ is unique -- when necessary, we denote it as $\omega_\ast(q)$. In the latter case, we denote $\omega'$ as $\omega'_\ast(q)$. Note in this case that either $\omega_\ast(q)=\omega'_\ast(q)2$ and $q=F_{\omega'_\ast(q)}(q_0)$, or symmetrically $\omega_\ast(q)=\omega'_\ast(q)3$ and $q=F_{\omega'_\ast(q)}(q'_0)$.

\subsection{Top balls}\label{top}
Consider a typical ball
$$B'_n:=B(q_1,2^{-n})=\bigcup\limits_{\omega\in\{0,1\}^n}\left(K_\omega\setminus F_\omega(\mathit{C})\right).$$
The following is a slightly modified analogue of \cite[Subsection 3.1]{HL}.

\begin{lemma}\label{exit-up}
$G_\mu^{B'_n}1(x)\asymp(2w_0)^nd_R^2\left(x,(B'_n)^c\right)$ for all $x\in B'_n$.
\end{lemma}

\begin{proof}
Set $\phi_n^x$ as the function attaining $R(x,(B'_n)^c)$.

To begin with, observe that
$$\phi_n^{q_1}=\begin{cases}
u_\downarrow\circ F_\omega^{-1}&\mbox{on}\quad K_\omega\quad\mbox{for each}\quad\omega\in\{0,1\}^n,\\
0&\mbox{elsewhere}.
\end{cases}$$
Therefore,
$$\mathcal{E}\left(\phi_n^{q_1}\right)=\sum\limits_{\omega\in\{0,1\}^n}s_\omega^{-1}\mathcal{E}(u\circ F_\omega)=2^n\cdot 2^n\mathcal{E}(u_\downarrow)=3\cdot 2^{2n}.$$
That is, $R(q_1,(B'_n)^c)\asymp 2^{-2n}$.

Next, we consider $y'_k=F_{0^n2^k}(q_1)$ (for all $k\ge 0$), which, due to reflection symmetry, stands for $F_{\omega\omega'}(q_1)$ with any $\omega\in\{0,1\}^n$ and $\omega'\in\{2,3\}^k$.

Denote for simplicity that
$$\phi_n^{y'_{k_0}}(y'_k)=b_k,\quad \phi_n^{y'_{k_0}}(F_{0^n2^k3}(q_1))=b'_k\quad\mbox{for all}\quad k\in\mathbb{N}.$$
Trivially $b_{k_0}=1$. Fix $k_0\ge 1$, then
$$\phi_n^{y'_{k_0}}=\begin{cases}b_0u_\downarrow\circ F_\omega^{-1}&\mbox{on}\quad K_\omega\quad\mbox{for each}\quad\omega\in\{0,1\}^n\setminus\{0^n\}\\
u_{b_{k+1},b_k,b_k}^-\circ F_{0^n2^k0}^{-1}&\mbox{on}\quad K_{0^n2^k0}\quad\mbox{for}\quad k=0,\dotsc,k_0-1\\
u_{b_k,b_k,b'_k}^-\circ F_{0^n2^k1}^{-1}&\mbox{on}\quad K_{0^n2^k1}\quad\mbox{for}\quad k=0,\dotsc,k_0-1\\
b'_ku_\downarrow\circ F_{0^n2^k3}^{-1}&\mbox{on}\quad K_{0^n2^k3}\quad\mbox{for}\quad k=0,\dotsc,k_0-1\\
u_\downarrow\circ F_{0^n2^{k_0}}^{-1}&\mbox{on}\quad K_{0^n2^{k_0}}\\
0&\mbox{elsewhere}\end{cases}$$
with $b_k,b'_k$ ($0\le k\le k_0-1$) optimizing the energy, that is, $b'_k=\frac{1}{4}b_k$ and
$$b_k=\frac{6\cdot 2^k-(6-7\cdot 2^{-n})4^{-k}}{6\cdot 2^{k_0}-(6-7\cdot 2^{-n})4^{-k_0}}\asymp\begin{cases}2^{-n-k_0}&\mbox{for}\quad k=0;\\
2^{k-k_0}&\mbox{for}\quad k=1,\dotsc,k_0.
\end{cases}$$
Repeating the procedure in \cite[Lemma 3.3]{HL}, we yield
$$R(y'_k,(B'_n)^c)\asymp 2^{-(n+k)},\quad\phi_n^{y'_k}(q_1)\asymp 2^{-(n+k)},\quad\phi_n^{y'_k}(y'_{k'})\asymp 2^{k'-k}\quad\mbox{for}\quad 0\le k'\le k.$$

Based on these estimates, we can easily pass to general $x\in B'_n$ in the same way as in \cite[Proposition 3.5 and Theorem 3.8]{HL}.
\end{proof}

\subsection{Unbalanced combination of cells}
\begin{figure}[htbp]
	\centering
	\begin{minipage}[t]{0.45\textwidth}
		\centering
		\includegraphics[scale=0.8]{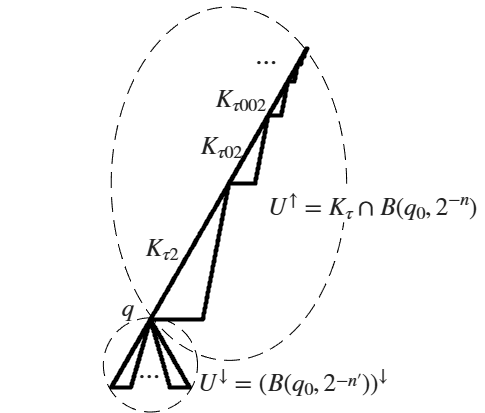}
		\caption{$U_{\tau,n'}^q$ of type a)}\label{U-a}
	\end{minipage}
	\begin{minipage}[t]{0.45\textwidth}
		\centering
		\includegraphics[scale=0.8]{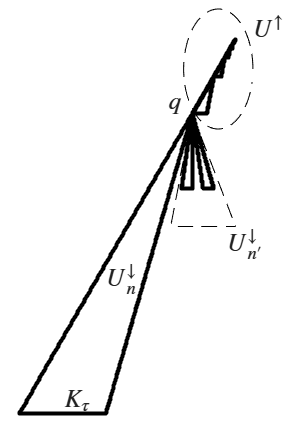}
		\caption{$U_{\tau,n'}^q$ of type b)}\label{U-b}
	\end{minipage}
\end{figure}

Here we focus on a new class of regions. Assume $q\in V_n^\uparrow$ with some $n\in\mathbb{N}$. For any $n'>n$ and $\tau\in\mathbb{S}_n$ satisfying $q\in K_\tau$, define
$$U_{\tau,n'}^q:=B\left(q,2^{-n'}\right)\cup\left(K_\tau\cap B(q,2^{-n})\right).$$
Such regions can be classified into two types (see Figures \ref{U-a} and \ref{U-b} respectively):
\begin{itemize}
\item[a)] $K_\tau$ lies above $q$. Here $U_{\tau,n'}^q$ contains $2^{n'-n_\ast(q)}$ almost disjoint $n'$-ordered cells below $q$ and a part of the $n$-ordered cell $K_\tau$. For simplicity, we set
$$U^\uparrow=K_\tau\cap B(q,2^{-n})\quad\mbox{and}\quad U^\downarrow=\cup\big\{K_\tau:\ \tau\in\mathbb{S}_{n'},\ q=F_\tau(q_1)\big\}.$$
Again we have $U^\uparrow\cup U^\downarrow=U_{\tau,n'}^q$ and $U^\uparrow\cap U^\downarrow=\{q\}$.
\item[b)] $K_\tau$ lies below $q$. Here $U_{\tau,n'}^q$ contains $K_\tau\setminus F_\tau(\mathit{C})$, $2^{n'-n_\ast(q)}-2^{n-n_\ast(q)}$ almost disjoint $n'$-ordered cells below $q$, and (if $n_\ast(q)>0$) a part of one $n'$-ordered cell above $q$. We denote these three parts respectively by $U^\downarrow_n$, $U^\downarrow_{n'}$ and $U^\uparrow$.
\end{itemize}

We consider case a) first, beginning by the basic location:

\begin{lemma}\label{G2}
For $U=U_{02^{n-1},n'}^{q_0}$,
$$G_\mu^U1(x)\asymp\begin{cases}
\left(2^nd_R(x,U^c)\right)^2\left(2^{-n'}w_0^{n'}+w_2^n\left(R(x,q_0)+2^{-2n'}\right)\right),&\mbox{for}\quad x\in U^\uparrow;\\
\left(2^{n'}d_R(x,U^c)\right)^2\left(2^{-n'}w_0^{n'}+2^{-2n'}w_2^n\right),&\mbox{for}\quad x\in U^\downarrow.
\end{cases}$$
\end{lemma}

\begin{proof}
Here $\tau=02^{n-1}$. Denote the function attaining $R(x,U)$ as $\psi_U^x$, and let
$$x_{m,k}=F_{02^{n-1}0^m23^k}(q_1)\quad\mbox{and}\quad y_k:=F_{20^{n'-1}2^k}(q_1).$$

For $x=q_0$, clearly
$$\psi_U^{q_0}=\begin{cases}
u_\downarrow\circ F_{2\omega}^{-1}&\mbox{on}\quad K_{2\omega}\quad\mbox{for each}\quad\omega\in\{0,1\}^{n'-1},\\
u_\uparrow\circ F_\tau^{-1}&\mbox{on}\quad K_\tau,\\
0&\mbox{elsewhere}.
\end{cases}$$
Thus
$$\mathcal{E}(\psi_U^{q_0})=\sum\limits_{\omega\in\{0,1\}^{n'-1}}s_{2\omega}^{-1}\mathcal{E}(u_\downarrow)+s_\tau^{-1}\mathcal{E}(u_\uparrow)=3\left(2^{2n'-1}+2^{n-1}\right)\asymp 2^{2n'}$$
and, using (\ref{int-psi}),
$$\int_{U_{\tau,n'}^{q_0}}\psi_U^{q_0}d\mu\asymp\mu\left(U_{\tau,n'}^{q_0}\right)=\sum\limits_{\omega\in\{0,1\}^{m-1}}w_{2\omega}+w_\tau
\asymp(2w_0)^{n'}+w_2^n.$$

For $x=x_{m_0,k_0}$, assume
$$\psi_U^{m_0,k_0}(q_0)=a_{-1,0},\quad \psi_U^{m_0,k_0}(F_{02^{n-1}0^m23^k}(q_1))=a_{m,k},\quad \psi_U^{m_0,k_0}(F_{02^{n-1}0^m23^k2}(q_1))=a'_{m,k}.$$
Then the construction of $\psi_U^{m_0,k_0}$ modifies that of $\psi_n^{x_{m_0,k_0}}$ only by $\omega\in\{0,1\}^{n'-1}$ replacing $\omega\in\{0,1\}^{n-1}$ on the first line. Accordingly, to optimize its energy, we modify the proof of \cite[Lemma 3.2]{HL} only by the first term
$$\frac{\partial\mathcal{E}}{\partial a_{-1,0}}=3\cdot 2^{2n'}a_{-1,0}+2\cdot 2^{n+1}(a_{-1,0}-a_{0,0})=0.$$
Thus in parallel,
\begin{align*}
a_{m,0}=&\ \frac{24\left(2^{2n'-n}+1\right)2^m-\left(3\cdot 2^{2n'-n}-4\right)4^{-m}}{12\left(2^{2n'-n}+1\right)\cdot 2^{m_0}\left(2^{k_0}+4^{-k_0}\right)-\left(3\cdot2^{2n'-n}-4\right)\cdot 4^{-m_0-k_0}}\\
a_{m_0,k}=&\ \frac{12\left(2^{2n'-n}+1\right)2^{m_0}\left(2^k+4^{-k}\right)-\left(3\cdot2^{2n'-n}-4\right)4^{-m_0-k}}{12\left(2^{2n'-n}+1\right)\cdot 2^{m_0}\left(2^{k_0}+4^{-k_0}\right)-\left(3\cdot2^{2n'-n}-4\right)\cdot 4^{-m_0-k_0}}.
\end{align*}
Hence for every $k_0\in\mathbb{N}$, for all $m=0,\dotsc,m_0$ and $k=1,\dotsc,k_0$,
$$a_{-1,0}\asymp 2^{n-2n'-m_0-k_0},\quad a_{m,0}\asymp a_{m,0}-a_{m-1,0}\asymp 2^{m-m_0-k_0}\quad\mbox{and}\quad a_{m_0,k}\asymp 2^{k-k_0}.$$
Consequently,
\begin{align*}
\mathcal{E}\left(\psi_U^{m_0,k_0}\right)\asymp&\ 2^{2n'}\cdot 4^{n-2n'-m_0-k_0}+2^{n+m_0+1}\cdot 4^{-k_0}+\sum\limits_{m=0}^{m_0}2^{n+m+1}\cdot 4^{m-m_0-k_0}\\
&\ +2^{n+m_0+1}\cdot 4^{-k_0}+2^{n+m_0+k_0+1}+\sum\limits_{k=1}^{k_0-1}2^{n+m_0+k+2}\cdot 4^{k-k_0}\\
\asymp&\ 2^{n+m_0+k_0}
\end{align*}
and (using (\ref{int-psi}) and the fact $0<2w_0=1-2w_2<1$)
\begin{align*}
\int_U\psi_U^{m_0,k_0}d\mu\asymp&\ a_{-1,0}\sum\limits_{\omega\in\{0,1\}^{n'-1}}\mu(K_{2\omega})+\sum\limits_{m=0}^{m_0-1}a_{m,0}w_0^mw_2^n+a_{m_0,0}w_0^{m_0}w_2^n\\
&\ +\sum\limits_{k=0}^{k_0-1}a_{m_0,k}w_0^{m_0}w_2^{n+k}+a_{m_0,k_0\vee 1}w_0^{m_0}w_2^{n+k_0}\\
\asymp&\ 2^{n-2n'-m_0-k_0}\cdot 2^{n'}w_0^{n'}+\sum\limits_{m=0}^{m_0}2^{m-m_0-k_0}w_0^mw_2^n+\sum\limits_{k=0}^{k_0}2^{k-k_0}w_0^{m_0}w_2^{n+k}\\
\asymp&\ 2^{-n'-m_0-k_0}\left(2^nw_0^{n'}+2^{n'}w_2^n\right).
\end{align*}

For $x=y_{k_0}$, let
$$\psi_U^{y_{k_0}}(y_k)=b_k,\quad \psi_U^{y_{k_0}}(F_{20^{n'-1}2^k3}(q_1))=b'_k\quad\mbox{for all}\quad k\in\mathbb{N}.$$
Then
$$\psi_U^{y_{k_0}}=\begin{cases}b_0u_\uparrow\circ F_\tau^{-1}&\mbox{on}\quad K_\tau\cap B_n\\
b_0u_\downarrow\circ F_{2\omega}^{-1}&\mbox{on}\quad K_{2\omega}\quad\mbox{for each}\quad\omega\in\{0,1\}^{n'-1}\setminus\{0^{n'-1}\}\\
u_{b_{k+1},b_k,b_k}^-\circ F_{20^{n'-1}2^k0}^{-1}&\mbox{on}\quad K_{20^{n'-1}2^k0}\quad\mbox{for}\quad k=0,\dotsc,k_0-1\\
u_{b_k,b_k,b'_k}^-\circ F_{20^{n'-1}2^k1}^{-1}&\mbox{on}\quad K_{20^{n'-1}2^k1}\quad\mbox{for}\quad k=0,\dotsc,k_0-1\\
b'_ku_\downarrow\circ F_{20^{n'-1}2^k3}^{-1}&\mbox{on}\quad K_{20^{n'-1}2^k3}\quad\mbox{for}\quad k=0,\dotsc,k_0-1\\
u_\downarrow\circ F_{20^{n'-1}2^{k_0}}^{-1}&\mbox{on}\quad K_{20^{n'-1}2^{k_0}}\\
0&\mbox{elsewhere}
\end{cases}$$
Optimizing its energy
\begin{align*}
\mathcal{E}\left(\psi_U^{y_{k_0}}\right)=&\ 2^n\cdot\frac{3}{2}b_0^2+(2^{n'-1}-1)2^{n'}\cdot 3b_0^2+2^{n'+k_0}\cdot 3\\
&\ +\sum\limits_{k=0}^{k_0-1}2^{n'+k+1}\left\{(b_k-b_{k+1})^2+(b_k-b'_k)^2+3(b'_k)^2\right\},
\end{align*}
we obtain $b'_k=\frac{1}{4}b_k$ and
$$b_k=\frac{3(1+2^{n-2n'})2^k-(3+3\cdot 2^{n-2n'}-7\cdot 2^{-n'})4^{-k}}{3(1+2^{n-2n'})2^{k_0}-(3+3\cdot 2^{n-2n'}-7\cdot 2^{-n'})4^{-k_0}}\asymp\begin{cases}2^{-n'-k_0}&(k=0);\\
2^{k-k_0}&(k=1,\dotsc,k_0).
\end{cases}$$
Therefore,
$$\mathcal{E}\left(\psi_U^{y_{k_0}}\right)\asymp \left(2^n+2^{2n'}\right)b_0^2+2^{n'+k_0}+\sum\limits_{k=0}^{k_0-1}2^{n'+k}b_k^2\asymp 2^{n'+k_0}$$
and
$$\int_U\psi_U^{y_{k_0}}d\mu\asymp 2^{-k_0}\left(2^{-n'}w_2^n+w_0^{n'}\right).$$

Treatment on $\psi_U^x$ for other points $x\in U$ repeats what was done in \cite[Remark 3.4, Proposition 3.5 and Theorem 3.8]{HL}.
\end{proof}

\begin{corollary}\label{G2+}
For every $U^q_{\tau,n'}$ with $q\in V_\ast$, $n'>n\ge n_\ast(q)$ and $\tau\in\mathbb{S}_n$ satisfying $q=F_\tau(q_2)$ or $F_\tau(q_3)$, then
$$G_\mu^U1\asymp\begin{cases}
w_{\omega_\ast(q)}\left(2^nd_R(x,U^c)\right)^2\left(2^{-n'}w_0^{n'-n_\ast(q)}+w_2^{n-n_\ast(q)}\left(R(x,q_0)+2^{-2n'}\right)\right)&\mbox{on}\quad U^\uparrow;\\
w_{\omega_\ast(q)}\left(2^{n'}d_R(x,U^c)\right)^2\left(2^{-n'}w_0^{n'-n_\ast(q)}+2^{-2n'}w_2^{n-n_\ast(q)}\right)&\mbox{on}\quad U^\downarrow.
\end{cases}$$
\end{corollary}

Since $R(x,q)\le 2^{-n}$, we have in particular that
\begin{equation}\label{U-up}
G_\mu^U1\lesssim w_{\omega_\ast(q)}\left(2^{-n'}w_0^{n'-n_\ast(q)}+2^{-n}w_2^{n-n_\ast(q)}\right).
\end{equation}

\begin{proof}
If $n_\ast(q)=0$, then $q=q_1$, for which the desired result is already clear by Lemma \ref{exit-up}. Otherwise, say $\omega_\ast(q)=\omega'_\ast2$, then $q=F_{\omega'_\ast}(q_0)$ and $\tau=\omega'_\ast\tau'$ with $\tau'\in\mathbb{S}_{n+1-n_\ast(q)}$. Thus by self-similarity,
$$\psi_{U^q_{\tau,n'}}^x=\psi_{F_{\omega'_\ast}^{-1}\left(U^q_{\tau,n'}\right)}^{F_{\omega'_\ast}^{-1}(x)}\circ F_{\omega'_\ast}=\psi_{U'}^{F_{\omega'_\ast}^{-1}(x)}\circ F_{\omega'_\ast},\quad\mbox{where}\quad U':=U_{\tau',n'+1-n_\ast(q)}^{q_0}.$$
Consequently, using Lemma \ref{G2}, we obtain for $q=F_\tau(q_2)$ that
\begin{align*}
G_\mu^U1(x)=&\ \frac{1}{\mathcal{E}\left(\psi_U^x\right)}\int_U\psi_U^xd\mu=\frac{1}{s_{\omega'_\ast}^{-1}\mathcal{E}\left(\psi_{U'}^{F_{\omega'_\ast}^{-1}(x)}\right)}\cdot w_{\omega'_\ast}\int\psi_{U'}^{F_{\omega'_\ast}^{-1}(x)}d\mu\\
=&\ 2^{-n_\ast(q)}w_{\omega'_\ast}G_\mu^{U'}1\left(F_{\omega'_\ast}^{-1}(x)\right)\\
\lesssim&\ 2^{-n_\ast(q)}w_{\omega'_\ast}\left(2^{-(n+1-n_\ast(q))}w_2^{n+1-n_\ast(q)}+2^{-(n'+1-n_\ast(q))}w_0^{n'+1-n_\ast(q)}\right)\\
\asymp&\ w_{\omega'_\ast}\left(2^{-n}w_2^{n-n_\ast(q)}+2^{-n'}w_0^{n'-n_\ast(q)}\right).
\end{align*}
The part ``$\gtrsim$'' is easy from the same construction.
\end{proof}

Using the same idea, we can easily show that

\begin{proposition}\label{G1}
For every $q\in V_n^\uparrow$, for all $x\in B:=B(q,2^{-n})$,
$$G_\mu^B1(x)\asymp\begin{cases}
2^nw_{\omega_\ast(q)}d_R^2(x,B^c)\left(w_0^{n-n_\ast(q)}+w_2^{n-n_\ast(q)}\left(R(x,q)+2^{-2(n-n_\ast(q))}\right)\right)&\mbox{on}\quad B^\uparrow;\\
2^nw_{\omega_\ast(q)}d_R^2(x,B^c)\left(w_0^{n-n_\ast(q)}+2^{-2n}w_2^{n-n_\ast(q)}\right)&\mbox{on}\quad B^\downarrow.
\end{cases}$$
\end{proposition}

This can also be considered as the case $n=n'$ in Corollary \ref{G2+}.\qed

Now we consider case b):

\begin{proposition}\label{G3}
For every $U^q_{\tau,n'}$ with $q\in V_\ast$, $n'>n\ge n_\ast(q)$ and $\tau\in\mathbb{S}_n$ satisfying $q=F_\tau(q_1)$, then
$$G_\mu^U1\asymp\begin{cases}
w_{\omega_\ast(q)}\left(2^nd_R(x,U^c)\right)^2\left(2^{-n}w_0^{n-n_\ast(q)}+2^{-2n'}w_2^{n'-n_\ast(q)}\right)&\mbox{on }\ U^\downarrow_n;\\
\begin{split}
w_{\omega_\ast(q)}\left(2^{n'}d_R(x,U^c)\right)^2\left\{2^{-2n'}w_0^{n-n_\ast(q)}+2^{n-2n'}w_0^{n'-n_\ast(q)}\right.\\
\left.\quad\quad+w_2^{n'-n_\ast(q)}\left(R(x,q_0)+2^{-2n'}\right)\right\}
\end{split}&\mbox{on }\ U^\uparrow;\\
w_{\omega_\ast(q)}\left(2^{n'}d_R(x,U^c)\right)^2\left(2^{-n'}w_0^{n'-n_\ast(q)}+2^{-2n'}\left(w_0^{n-n_\ast(q)}+w_2^{n'-n_\ast(q)}\right)\right)&\mbox{on }\ U^\downarrow_{n'}.
\end{cases}$$
\end{proposition}

Note that in comparison to (\ref{U-up}), we have
\begin{equation}\label{U-down}
G_\mu^U1\lesssim w_{\omega_\ast(q)}\left(2^{-n}w_0^{n-n_\ast(q)}+2^{-n'}w_2^{n'-n_\ast(q)}\right).
\end{equation}

\begin{proof}
We only provide calculus on $U=U_{\tau,n'}^{q_0}$ with $\tau=20^{n-1}$. For other locations, if $n_\ast(q)>0$, we use the same method as in Corollary \ref{G2+}; while for $q=q_1$, we modify the proof in the same way as in Subsection \ref{top}.

Denote also
$$\Lambda_{n'}=\left\{\omega=i_2\cdots i_{n'}\in\{0,1\}^{n'-1}:\ i_2\cdots i_n\ne 0^{n-1}\right\},$$
and set
$$x_{m,k}=F_{02^{n'-1}0^m23^k}(q_1),\quad y_k:=F_{20^{n-1}2^k}(q_1)\quad\mbox{and}\quad z_k:=F_{21^{n'-1}2^k}(q_1).$$

For $x=q_0$, clearly
$$\psi_U^{q_0}=\begin{cases}
u_\downarrow\circ F_\tau^{-1}&\mbox{on}\quad K_\tau,\\
u_\downarrow\circ F_{2\omega}^{-1}&\mbox{on}\quad K_{2\omega}\quad\mbox{for each}\quad\omega\in\Lambda_{n'},\\
u_\uparrow\circ F_{02^{n'-1}}^{-1}&\mbox{on}\quad K_{02^{n'-1}},\\
0&\mbox{elsewhere}.
\end{cases}$$
Thus
$$\mathcal{E}(\psi_U^{q_0})=\frac{3}{2}\cdot 2^n+\frac{3}{2}\cdot 2^{n'}\left(2^{n'-1}-2^{n'-n}\right)+3\cdot 2^{n'}\asymp 2^{2n'}$$
and
\begin{align*}
\int_{U_{\tau,n'}^{q_0}}\psi_U^{q_0}d\mu\asymp&\ \mu\left(U_{\tau,n'}^{q_0}\right)=w_2w_0^{n-1}+\left(2^{n'-1}-2^{n'-n}\right)w_2w_0^{n'-1}+w_0w_2^{n'-1}\\
\asymp&\ w_0^n+(2w_0)^{n'}+w_2^{n'}.
\end{align*}

For $x=x_{m_0,k_0}$, the construction is similar as before. The first restriction is modified to
$$\frac{\partial\mathcal{E}}{\partial a_{-1,0}}=3\cdot 2^{n'}\left(2^{n'-1}-2^{n'-n}\right)a_{-1,0}+2^{n'+1}(a_{-1,0}-a_{0,0})=0,$$
yielding values of $\psi_U^{m_0,k_0}$ at $x_{m,0}$ and $x_{m_0,k}$ by
\begin{align*}
a_{m,0}=&\ \frac{12\left(2\left(2^{n'-1}-2^{n'-n}\right)+1\right)2^m-\left(3\left(2^{n'-1}-2^{n'-n}\right)-2\right)4^{-m}}{6\left(2\left(2^{n'-1}-2^{n'-n}\right)+1\right)\cdot 2^{m_0}\left(2^{k_0}+4^{-k_0}\right)-\left(3\left(2^{n'-1}-2^{n'-n}\right)-2\right)\cdot 4^{-m_0-k_0}};\\
a_{m_0,k}=&\ \frac{6\left(2\left(2^{n'-1}-2^{n'-n}\right)+1\right)2^{m_0}\left(2^k+4^{-k}\right)-\left(3\left(2^{n'-1}-2^{n'-n}\right)-2\right)4^{-m_0-k}}{6\left(2\left(2^{n'-1}-2^{n'-n}\right)+1\right)\cdot 2^{m_0}\left(2^{k_0}+4^{-k_0}\right)-\left(3\left(2^{n'-1}-2^{n'-n}\right)-2\right)\cdot 4^{-m_0-k_0}}.
\end{align*}
In particular,
$$a_{-1,0}\asymp 2^{-n'-m_0-k_0},\quad a_{m,0}\asymp a_{m,0}-a_{m-1,0}\asymp 2^{m-m_0-k_0}\quad\mbox{and}\quad a_{m_0,k}\asymp 2^{k-k_0}.$$
Consequently,
$$\mathcal{E}\left(\psi_U^{m_0,k_0}\right)\asymp 2^{n'+m_0+k_0}$$
and (using the fact $0<2w_0=1-2w_2<1$)
\begin{align*}
\int_U\psi_U^{m_0,k_0}d\mu\asymp&\ 2^{-n'-m_0-k_0}\cdot\left(w_0^n+\left(2^{n'-1}-2^{n'-n}\right)w_0^{n'}\right)\\
&\ +\sum\limits_{m=0}^{m_0}2^{m-m_0-k_0}w_0^mw_2^{n'}+\sum\limits_{k=0}^{k_0}2^{k-k_0}w_0^{m_0}w_2^{n'+k}\\
\asymp&\ 2^{-(n'+m_0+k_0)}\left(w_0^n+2^nw_0^{n'}+2^{n'}w_2^{n'}\right).
\end{align*}

For $x=y_{k_0}$, there is
$$\psi_U^{y_{k_0}}=\begin{cases}b_0u_\uparrow\circ F_{02^{n'-1}}^{-1}&\mbox{on}\quad K_{02^{n'-1}}\cap B'_n\\
b_0u_\downarrow\circ F_{2\omega}^{-1}&\mbox{on}\quad K_{2\omega}\quad\mbox{for each}\quad\omega\in\Lambda_{n'}\\
u_{b_{k+1},b_k,b_k}^-\circ F_{20^{n-1}2^k0}^{-1}&\mbox{on}\quad K_{20^{n-1}2^k0}\quad\mbox{for}\quad k=0,\dotsc,k_0-1\\
u_{b_k,b_k,b'_k}^-\circ F_{20^{n-1}2^k1}^{-1}&\mbox{on}\quad K_{20^{n-1}2^k1}\quad\mbox{for}\quad k=0,\dotsc,k_0-1\\
b'_ku_\downarrow\circ F_{20^{n-1}2^k3}^{-1}&\mbox{on}\quad K_{20^{n-1}2^k3}\quad\mbox{for}\quad k=0,\dotsc,k_0-1\\
u_\downarrow\circ F_{20^{n-1}2^{k_0}}^{-1}&\mbox{on}\quad K_{20^{n-1}2^{k_0}}\\
0&\mbox{elsewhere}
\end{cases}$$
with
$$b_k=\psi_U^{y_{k_0}}(y_k)\quad\mbox{and}\quad b'_k=\psi_U^{y_{k_0}}(F_{20^{n-1}2^k3}(q_1))=\frac{1}{4}b_k.$$
Similar optimization as before yields
$$b_k=\frac{\left(3\left(2^{n'}-2^{n'-n+1}+1\right)+6\cdot 2^{n-n'}\right)2^k-\left(3\left(2^{n'}-2^{n'-n+1}+1\right)-2^{n-n'}\right)4^{-k}}{\left(3\left(2^{n'}-2^{n'-n+1}+1\right)+6\cdot 2^{n-n'}\right)2^{k_0}-\left(3\left(2^{n'}-2^{n'-n+1}+1\right)-2^{n-n'}\right)4^{-k_0}}.$$
That is,
$$b_0\asymp 2^{n-2n'-k_0}\quad\mbox{and}\quad b_k\asymp b_k-b_{k-1}\asymp 2^{k-k_0}\quad(k=1,\dotsc,k_0).$$
It follows that
$$\mathcal{E}\left(\psi_U^{y_{k_0}}\right)\asymp 2^{n+k_0}\quad\mbox{and}\quad\int_U\psi_U^{y_{k_0}}\asymp 2^{-n-k_0}\left(2^{2(n-n')}w_2^{n'}+2^{n-n'}w_0^{n'}+2^nw_0^n\right).$$

Moreover, for $x=z_k$, there is
$$\psi_U^{z_{k_0}}=\begin{cases}b_0u_\uparrow\circ F_{02^{n'-1}}^{-1}&\mbox{on}\quad K_{02^{n'-1}}\cap B'_n\\
b_0u_\downarrow\circ F_{2\omega}^{-1}&\mbox{on}\quad K_{2\omega}\quad\mbox{for each}\quad\omega\in\left\{0^{n-1}\right\}\cup\Lambda_{n'}\\
u_{b_{k+1},b_k,b_k}^-\circ F_{20^{n'-1}3^k0}^{-1}&\mbox{on}\quad K_{20^{n'-1}3^k0}\quad\mbox{for}\quad k=0,\dotsc,k_0-1\\
u_{b_k,b_k,b'_k}^-\circ F_{20^{n'-1}3^k1}^{-1}&\mbox{on}\quad K_{20^{n'-1}3^k1}\quad\mbox{for}\quad k=0,\dotsc,k_0-1\\
b'_ku_\downarrow\circ F_{20^{n'-1}3^k2}^{-1}&\mbox{on}\quad K_{20^{n'-1}3^k2}\quad\mbox{for}\quad k=0,\dotsc,k_0-1\\
u_\downarrow\circ F_{20^{n'-1}3^{k_0}}^{-1}&\mbox{on}\quad K_{20^{n'-1}3^{k_0}}\\
0&\mbox{elsewhere}
\end{cases}$$
with
$$b_k=\psi_U^{y_{k_0}}(y_k)\quad\mbox{and}\quad b'_k=\psi_U^{y_{k_0}}(F_{20^{n'-1}3^k2}(q_1))=\frac{1}{4}b_k.$$
Optimizing its energy
\begin{align*}
\mathcal{E}\left(\psi_U^{z_{k_0}}\right)=&\ 2^{n'}\cdot\frac{3}{2}b_0^2+\left(2^n+\left(2^{n'-1}-2^{n'-n}-1\right)2^{n'}\right)\cdot 3b_0^2+2^{n'+k_0}\cdot 3\\
&+\sum\limits_{k=0}^{k_0-1}2^{n'+k+1}\left\{(b_k-b_{k+1})^2+\frac{3}{4}b_k^2\right\},
\end{align*}
we obtain
$$b_k=\frac{\left(3\Delta_{n,n'}+6\right)2^k-\left(3\Delta_{n,n'}-1\right)4^{-k}}{\left(3\Delta_{n,n'}+6\right)2^{k_0}-\left(3\Delta_{n,n'}-1\right)4^{-k_0}},\quad\mbox{where}\quad\Delta_{n,n'}=2^{n'}-2^{n'-n+1}+2^{n+1-n'}-1.$$
That is,
$$b_0\asymp 2^{-n'-k_0}\quad\mbox{and}\quad b_k\asymp b_k-b_{k-1}\asymp 2^{k-k_0}\quad(k=1,\dotsc,k_0).$$
It follows that
$$\mathcal{E}\left(\psi_U^{z_{k_0}}\right)\asymp 2^{n'+k_0}\quad\mbox{and}\quad\int_U\psi_U^{z_{k_0}}\asymp 2^{-n'-k_0}\left(2^{n'}w_0^{n'}+w_2^{n'}+w_0^n\right).$$

The desired estimates are then proved for typical points $q_0$, $x_{m,k}$, $y_k$ and $z_k$. The counterpart for general $x\in U$ follows similarly as in Lemma \ref{G2}.
\end{proof}

\section{Exit time for general balls in the case $w_2\le w_0$}\label{S5}

Throughout this section we assume $w_2\le w_0$. In this setting, we construct scales related to $(K,R,\mu,\mathcal{E},\mathcal{F})$.

We first replace $n$ in Section \ref{S4} by non-integer's. Note that it can not be obtained from easy interpolation technique like taking $2^{-n-1}<r\le 2^{-n}$, since we need lower estimate on $B(q,(1-\varepsilon_0)r)\setminus B(q,2^{-n-1})$ also.

\begin{proposition}\label{G'r-1}
For every $\varepsilon_0\in(0,1)$, there exists $\delta_0>0$ such that, for every ball $B=B(q,r)$ with $q\in V_n^\uparrow$ and $0<r\le 2^{-n}$,
$$\delta_01_{(1-\varepsilon_0)B}\lesssim\frac{G_\mu^B1}{w_{\omega_\ast(q)}w_0^{-n_\ast(q)}r^{\log_2w_0^{-1}+1}}\lesssim 1.$$
\end{proposition}

\begin{proof}
The desired result is provided by Proposition \ref{G1} if $r=2^{-(n+m)}$ for some $m\in\mathbb{N}$. Now we consider other $r$.

Take the unique $n'\in\mathbb{N}$ such that $2^{-n'-1}<r\le 2^{-n'}$. Clearly $n'\ge n$, and hence $q\in V_n^\uparrow\subset V_{n'}^\uparrow$. Set $B'=B(q,2^{-n'})$. Using the assumption $w_2\le w_0$ and the fact $R(x,q)<r\le 2^{-n'}$, we see directly that
$$G_\mu^{B'}1(x)\lesssim 2^{n'}w_{\omega_\ast(q)}d_R^2\left(x,(B')^c\right)w_0^{n'-n_\ast(q)}\le 2^{-n'}w_{\omega_\ast(q)}w_0^{n'-n_\ast(q)}.$$
From the selection of $n'$ we see $n'\asymp\log_2r^{-1}$. Therefore,
\begin{align*}
G_\mu^B1\le&\ G_\mu^{B'}1\lesssim w_{\omega_\ast(q)}w_0^{-n_\ast(q)}\left(\frac{w_0}{2}\right)^{n'}\\
\asymp&\ w_{\omega_\ast(q)}w_0^{-n_\ast(q)}\left(\frac{w_0}{2}\right)^{\log_2r^{-1}}=w_{\omega_\ast(q)}w_0^{-n_\ast(q)}r^{\log_2w_0^{-1}+1}.
\end{align*}
The upper estimate is then proved.

Now we turn to the lower side. For one hand, observe that
$$\phi:=\left(1-\sup_{B'\setminus B}\psi_{B'}^x\right)^{-1}\left(\psi_{B'}^x(y)-\sup_{B'\setminus B}\psi_{B'}^x\right)_+$$
is a test function for $R(x,B^c)$. From the proof of Lemma \ref{exit} we can easily check that for every $\varepsilon_0\in(0,1)$, there exists $\varepsilon'_0\in(0,1)$ (independent of $B'$) such that
$$\psi_{B'}^x(y)\le\varepsilon'_0\quad\mbox{for all}\quad x\in(1-\varepsilon_0)B'\quad\mbox{and}\quad y\in B'\setminus\varepsilon_0B'.$$
Hence
\begin{align}
R(x,B^c)\ge&\ \frac{1}{\mathcal{E}(\phi)}=\frac{\left(1-\sup_{B'\setminus B}\psi_{B'}^x\right)^2}{\mathcal{E}\left(\left(\psi_{B'}^x-\sup_{B'\setminus B}\psi_{B'}^x\right)_+\right)}\ge\frac{\left(1-\varepsilon'_0\right)^2}{\mathcal{E}(\psi_{B'}^x)}\notag\\
\gtrsim&\ d_R(x,(B')^c)\ge d_R(x,B^c)\ge\varepsilon_0r.\label{R'low}
\end{align}
Recall that $R$ is geodesic on $K$. Let $\overline{xy}$ be the unique path from $x$ to $y$. Then by the maximum principle, for every $z\in\overline{xq}$,
\begin{equation}\label{cmprs}
\psi_B^x=\psi_B^x(z)\psi_B^z
\end{equation}
on the connected component of $B\setminus\{z\}$ containing $q$.

Now we iterate (\ref{cmprs}) to obtain lower estimate of $\psi_B^x$. Before that we need some preparations. Let $\omega=\omega_\ast(q)=\emptyset$ if $q=q_1$, while take $\omega=\omega'_\ast(q)$ if $q\ne q_1$ (say $\omega_\ast(q)=\omega 2$). Assume the dyadic representations of $2^{n'}r$ and $2^{n'}(1-\frac{\varepsilon_0}{2})r$ are respectively $0.a_1a_2\cdots$ and $0.a'_1a'_2\cdots$. Note that $a_1=1$. Let $n''$ be the first integer such that $a_{n''}\ne a'_{n''}$ (clearly $a_{n''}=1$, $a'_{n''}=0$), and $m_0$ be the last integer such that $1-2^{-m_0}<2^{n'}r$. It is easy to see
\begin{equation}\label{n''+}
n'+n''\le 1+\log_2\left((\varepsilon_0r)^{-1}\right).
\end{equation}

We only consider the case where $x$ lies on a geodesic from $q$ to $\partial B$ (then, for general $x\in B$, the desired result holds for $\bar{x}$, and the counterpart for $x$ follows by repeating the last step in the proof of \cite[Proposition 3.5]{HL}). Further, by symmetry, we assume without loss of generality that $x$ lies between ``typical'' points of $B(q,2^{-n})$. Under this assumption, let us consider $B_{0,k}$, the farthest ball -- centred at a $V_{n+k}^\uparrow$ point on typical segments, of radius $2^{-n-k}$ -- contained in $B(q,r)$.

To be exact, set $q_{0,k}=q$ whenever $k\le 0$. For every $x\in B^\downarrow$ and $1\le k\le n''$, set $i_{0,k}=0$ if $a_k=0$, while $i_{0,k}=2$ if $a_k=1$. Take $q_{0,k}=F_{\omega_\ast(q)0^{n'-n_\ast(q)}i_{0,1}\cdots i_{0,k}}(q_1)$ if $a_{k+1}=1$; while take $q_{0,k}=F_{\omega_\ast(q)0^{n'-n_\ast(q)}i_{0,1}\cdots i_{0,k'-1}02^{k-k'+1}}(q_1)$ if $a_{k+1}=0$ and $a_{k'}=1$, $a_{k''}=0$ for all $k''\in(k',k]$ with some $k'\le k$. Similarly, for every $x\in B^\uparrow$, fix the unique $m\ge 0$ such that $x\in\overline{x_{m,k_0-1}x_{m,k_0}}$ for some $k_0\ge 0$ (here $q$ is considered as $x_{0,-1}$). For all $k\le m$, set $q_{0,k}=x_{k,0}=F_{\omega02^{n'-n_\ast(q)}0^m2}(q_1)$; while for all $m<k\le k_0$, set $i_{0,k}=1$ if $a_k=0$; $i_{0,k}=3$ if $a_k=1$. Take $q_{0,k}=F_{\omega02^{n'-n_\ast(q)}0^m2i_{0,m+1}\cdots i_{0,k}}(q_1)$ if $a_{k+1}=1$; while take $q_{0,k}=F_{\omega02^{n'-n_\ast(q)}0^m2i_{0,m+1}\cdots i_{0,k'-1}13^{k-k'+1}}(q_1)$ if $a_{k+1}=0$ and $a_{k'}=1$, $a_{k''}=0$ for all $k''\in(k',k]$ with some $k'\le k$.

\begin{figure}[htbp]
\centering\includegraphics[width=0.4\textwidth]{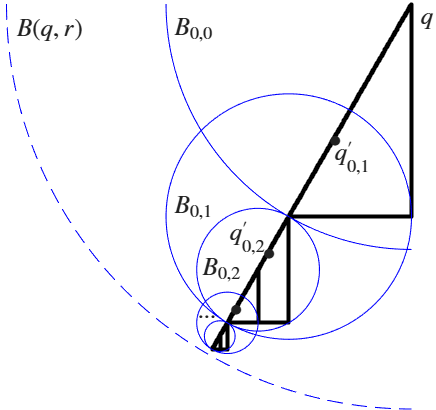}
\caption{Balls $B_{0,k}$ and the intermediate points $q'_{0,k}$ for $r=\frac{4}{5}\cdot 2^{-n}$.}
\label{kLay}
\end{figure}

In all cases above, define $B_{0,k}=B(q_{0,k},2^{-n'-k-1})$. See Figure \ref{kLay} for an illustration of these balls. Take integer $0<k'\le n''$ such that $x\in\overline{q_{0,k'-1}q_{0,k'}}$. By geometry, if $k'<n''$, then $x$ is in either $\frac{9}{10}B_{0,k'-1}$ or $\frac{9}{10}B_{0,k'}$ (say the former); while if $k'=n''$, then $x\in(1-\varepsilon_0)B_{0,k'-1}$ by the inequalities on radii.

Further, let $q'_{0,k}$ be the point on $\overline{q_{0,k-1}q_{0,k}}$ such that $R(q'_{0,k},q_{0,k})=\frac{1}{3}R(q_{0,k-1},q_{0,k})$. In particular, $q'_{0,k}\in\frac{9}{10}B_{0,k-1}\cap\frac{9}{10}B_{0,k}$. Then
$$\psi_{B_{0,k'-1}}^x(q'_{0,k'-1})\ge\varepsilon_1\quad\mbox{and}\quad\psi_{B_{0,k}}^{q'_{0,k}}(q'_{0,k-1})\ge\varepsilon_1\quad\mbox{for all}\quad 1\le k\le n''.$$
with a constant $\varepsilon_1$ (independent of $B,x,k,k'$). It follows for all $y\in B(q,2^{-n'-2})$ that
\begin{align*}
\psi_B^x(y)\ge&\ \psi_B^x(q'_{0,k'-2})\ \psi_B^{q'_{0,k'-2}}(y)\ge\psi_B^x(q'_{0,k'-2})\ \psi_B^{q'_{0,k'-2}}(q'_{0,k'-3})\ \psi_B^{q'_{0,k'-3}}(y)\ge\cdots\\
\ge&\ \psi_B^x(q'_{0,k'-2})\ \psi_B^{q'_{0,k'-2}}(q'_{0,k'-3})\cdots\psi_B^{q'_{0,1}}(q'_{0,0})\ \psi_B^{q'_{0,0}}(y)\\
\ge&\ \psi_{B_{0,k'-1}}^x(q'_{0,k'-2})\ \psi_{B_{0,k'-2}}^{q'_{0,k'-2}}(q'_{0,k'-3})\cdots\psi_{B_{0,1}}^{q'_{0,1}}(q'_{0,0})\ \psi_{B_{0,0}}^{q'_{0,0}}(y)\\
\ge&\ \varepsilon_1^{k'}\ge\varepsilon_1^{n''}\ge\varepsilon_1^{1+\log_2\left((\varepsilon_0r)^{-1}\right)-n'}\ge\varepsilon_2\varepsilon_1^{\log_2(\varepsilon_0^{-1})}
\end{align*}
with a constant $\varepsilon_2$ depending only on $\varepsilon_0,\varepsilon_1$. Here we use (\ref{n''+}) for the last but one inequality. Combining (\ref{R'low}), we obtain (with $\delta_0=\varepsilon_0\varepsilon_2\varepsilon_1^{\log_2(\varepsilon_0^{-1})}$)
$$G_\mu^B1(x)\ge R(x,B^c)\int_{B(q,2^{-n'-2})}\psi_B^xd\mu\ge\varepsilon_0r\cdot\varepsilon_1^{k'}\mu(B(q,2^{-n'-2}))\gtrsim\delta_0w_{\omega_\ast(q)}w_0^{-n_\ast(q)}r^{\log_2w_0^{-1}+1},$$
proving the desired lower estimate exactly.
\end{proof}

The same can be easily generalized to regions interpolating $U_{\tau,n'}^q$:

\begin{proposition}\label{G'r-2}
For every $\varepsilon_0\in(0,1)$, there exists $\delta_0>0$ such that, for every $q\in V_n^\uparrow$ and $0<r\le 2^{-n}$, if $q\in F_\tau(V_0)$ with some $\tau\in\mathbb{S}_n$, then, given every $n'>\log_2(\varepsilon_0r)^{-1}$, for $U=B(q,2^{-n'})\cup(B(q,r)\cap K_\tau)$,
$$\delta_01_{U_{\varepsilon_0}}\lesssim\frac{G_\mu^U1}{w_{\omega_\ast(q)}\left(2^{-n'}w_0^{n'-n_\ast(q)}+w_2^{-n_\ast(q)}r^{\log_2w_2^{-1}+1}\right)}\lesssim 1$$
if $q\in\{F_\tau(q_2),F_\tau(q_3)\}$, while
$$\delta_01_{U_{\varepsilon_0}}\lesssim\frac{G_\mu^U1}{w_{\omega_\ast(q)}\left(2^{-2n'}w_2^{n'-n_\ast(q)}+w_0^{-n_\ast(q)}r^{\log_2w_0^{-1}+1}\right)}\lesssim 1$$
if $q=F_\tau(q_1)$. Here $U_{\varepsilon_0}=\left\{x\in U:d_R(x,U^c)\ge\varepsilon_0r\right\}$.\qed
\end{proposition}

Now we are ready for general balls. Given arbitrary $x\in K$ and $r\in(0,1)$, let $n_\ast(x,r)$ be the least integer $n$ such that there exists $\omega_\ast(x,r)\in\mathbb{S}^n$ satisfying $q_\ast(x,r):=F_{\omega_\ast(x,r)}(q_1)\in B(x,r)$. Here are some basic observations:
\begin{itemize}
\item $n_\ast(x,r)=n_\ast(q_\ast(x,r))$ and $\omega_\ast(x,r)=\omega_\ast(q_\ast(x,r))$. In particular, either $n_\ast(x,r)=0$, or $i_{n_\ast(x,r)}\in\{2,3\}$;
\item $n_\ast(q,r)=n_\ast(q)$ and $\omega_\ast(q,r)=\omega_\ast(q)$ for all $r\le 2^{-n}$, if $q\in V_n^\uparrow$;
\item $B(x,r)\subset B\left(q_\ast(x,r),2^{-n_\ast(x,r)}\right)$. In particular, $\omega_\ast(x,r)$ and $q_\ast(x,r)$ are unique;
\item $n_\ast(q_1,r)=0$ and $n_\ast(q_2,r)=n_\ast(q_3,r)\asymp\log_2r^{-1}$.
\end{itemize}

Now we introduce the main theorem of this paper:

\begin{theorem}\label{EW}
Assume that $w_2\le w_0$. Then, $(\mathrm{E}_W)$ holds with $\sigma=\frac{1}{2}$, $\delta=\frac{1}{5}$ and
$$W(x,r)=\begin{cases}
w_{\omega_\ast(x,r)}w_0^{-n_\ast(x,r)}r^{\log_2w_0^{-1}+1},&\mbox{if}\quad x\in K_{\omega_\ast(x,r)};\\
w_{\omega_\ast(x,r)}\left(w_2^{-n_\ast(x,r)}r^{\log_2w_2^{-1}+1}+w_0^{-n_\ast(q)}\left\{r-R(x,q_\ast(x,r))\right\}^{\log_2w_0^{-1}+1}\right),&\mbox{if}\quad x\notin K_{\omega_\ast(x,r)}.
\end{cases}.$$
\end{theorem}

\begin{proof}
Fix $B=B(x,r)$ with $r\le 1$, and write $q_\ast=q_\ast(x,r)$, $n_\ast=n_\ast(x,r)$ for simplicity. We distinguish the following two cases:

Case 1: $q_\ast\in\frac{1}{5}B$. By applying Proposition \ref{G'r-1} to $B':=B\left(q_\ast,\frac{4}{5}r\right)$ we obtain
$$\inf\limits_{\frac{1}{5}B}G_\mu^B1\ge\inf\limits_{\frac{1}{2}B'}G_\mu^{B'}1\gtrsim w_{\omega_\ast(q_\ast)}w_0^{-n_\ast(q_\ast)}\left(\frac{4}{5}r\right)^{\log_2w_0^{-1}+1}\asymp w_{\omega_\ast(x,r)}w_0^{-n_\ast(x,r)}r^{\log_2w_0^{-1}+1};$$
and to $B'':=B\left(q_\ast,\frac{6}{5}r\wedge 2^{-n_\ast}\right)$ we obtain
$$\sup\limits_BG_\mu^B1\le\sup\limits_BG_\mu^{B''}1\lesssim w_{\omega_\ast(q_\ast)}w_0^{-n_\ast(q_\ast)}\left(\frac{6}{5}r\right)^{\log_2w_0^{-1}+1}\asymp w_{\omega_\ast(x,r)}w_0^{-n_\ast(x,r)}r^{\log_2w_0^{-1}+1}.$$
The desired result gets proved clearly.
	
Case 2: $q_\ast\in B\setminus\frac{1}{5}B$. Take the unique $n''\in\mathbb{N}$ such that
$$2^{-(n'+1)}<r-R(x,q_\ast)\le 2^{-n'}.$$
By geometry we can easily check that
$$\frac{1}{5}B\subset U_-\subset B\subset U_+,$$
where
\begin{align*}
U_-=&\ B\left(q_\ast,2^{-(n'+1)}\right)\cup\left(K_{\omega 0}\cap B\left(q_\ast,R(x,q_\ast)+\frac{r}{5}\right)\right);\\
U_+=&\ B\left(q_\ast,2^{-n'}\right)\cup\left(K_{\omega 0}\cap B(q_\ast,R(x,q_\ast)+r)\right).
\end{align*}
Note that $x\in K_{\omega_\ast(x,r)}$ corresponds to the case $q_\ast=F_{\omega_\ast(x,r)}(q_1)$, while $x\notin K_{\omega_\ast(x,r)}$ corresponds to the case $q_\ast=F_{\omega_\ast(x,r)}(q_2)$ or $F_{\omega_\ast(x,r)}(q_3)$. Therefore, by applying Proposition \ref{G'r-2} to $U_\pm$ we yield the desired estimates parallel to Case 1.
\end{proof}

However:

\begin{lemma}
$W$ is not a scaling function whenever $w_2<w_0$.
\end{lemma}

\begin{proof}
Let $y_n=F_{02^{n-1}}$ and $r_n=2^{-n}$. Clearly $q_\ast(y_n,r_n)=y_n$ and $q_\ast(y_n,2r_n)=q_0$. Therefore, by Theorem \ref{EW},
$$W(y_n,r_n)=w_{02^{n-1}}w_0^{-n}r_n^{\log_2w_0^{-1}+1}=w_0w_2^{n-1}2^{-n}$$
and
$$W(y_n,2r_n)=w_2\left(w_2^{-1}(2r_n)^{\log_2w_2^{-1}+1}+w_0^{-1}(2r_n-r_n)^{\log_2w_0^{-1}+1}\right)\asymp2^{-n}w_0^n.$$
It follows that
$$\frac{W(y_n,2r_n)}{W(y_n,r_n)}\asymp\left(\frac{w_0}{w_2}\right)^n\to\infty\quad\mbox{as}\quad n\to\infty,$$
contradicting the scaling property (\ref{scal}). Thus $W$ is not a scaling function.
\end{proof}

On the other hand, if $w_0=w_2=\frac{1}{4}$, then $W(x,r)=r^3$. This is clearly a scaling function, and coincides with the scale of a normal interval or a non-fractal tree. However, we still have $(\mathrm{wEH})\setminus(\mathrm{EHI})$ due to \cite[Theorems 4.1 and 4.2]{HL}. This shows that the known implication $(\mathrm{EHI})+(\mathrm{E})\Leftrightarrow(\mathrm{wEH})+(\mathrm{E})$ requires $(\mathrm{VD})$ essentially.

Now we summarize the important properties depending on the ratio $w_0/w_2$ (from \cite[Theorems 3.9, 4.1 and 4.2]{HL}and our Theorem 4.3, Lemma 4.4):
\begin{table}[htbp]
\centering
\caption{Different phenomenon as the ratio $w_0/w_2$ varies (note that $(\mathrm{EHI})$ always fails)}
\begin{tabular}{ccc}
\toprule
Weights&Harnack property&Exit time estimate\\
\midrule
$w_2<w_0$&$(\mathrm{wEH}_1)$ holds&$(\mathrm{E}_W)$ holds but $\frac{W(x,2r)}{W(x,r)}$ is unbounded\\
$w_2=w_0$&$(\mathrm{wEH}_1)$ holds&$(\mathrm{E}_W)$ holds with scaling function $W$\\
$w_0<w_2<2w_0$&$(\mathrm{wEH}_\gamma)$ with $\gamma<1$&$(\mathrm{E}_\bullet)$ fails since $\frac{\sup_BG_\mu^B1}{\inf_{\delta B}G_\mu^B1}$ is unbounded\\
$w_2\ge 2w_0$&$(\mathrm{wEH}_\bullet)$ fails&$(\mathrm{E}_\bullet)$ fails since $\frac{\sup_BG_\mu^B1}{\inf_{\delta B}G_\mu^B1}$ is unbounded\\
\bottomrule
\end{tabular}
\end{table}

\section{Capacity and heat kernel estimates}
\subsection{Capacity estimate}
Based on Theorem \ref{EW}, we can easily obtain survival and capacity estimates using (\ref{imply}), with respect to the singular scale $W$. Further, in the same way as how $R(x,B^c)$ was estimated, we can directly verify a much stronger (actually, the optimal) capacity estimate:

\begin{lemma}\label{cap+}
Given any $\lambda\in(0,1)$, for every ball $B=B(x,r)$ with $r\in(0,1]$,
$$\mathrm{cap}_0(\lambda B,B)\asymp\frac{1}{\gamma_\lambda(x,r)},$$
where
$$\gamma_\lambda(x,r)=2^{n_\ast(x,r)}r\left\{d_R\big(q_\ast(x,r),B^c\big)\wedge\big((1-\lambda)r\big)\right\}+d_R\big(q_\ast(x,r),\lambda B\big).$$
\end{lemma}

Note that we assume nothing on $w_0$ and $w_2$ in Lemma \ref{cap+}. Actually, this result concerns the resistance form only, without any reference measure involved.

Here are some observations on $\gamma_\lambda$:
\begin{itemize}
\item if $q_\ast(x,r)\in\lambda B$, then $\gamma_\lambda(x,r)=(1-\lambda)2^{n_\ast(x,r)}r^2$. Further, fixing any $\lambda'\in(\lambda,1)$, $\gamma_\lambda$ is comparable to the same term when $q_\ast(x,r)\in\lambda'B$;
\item $(1-\lambda)r^2\le\gamma_\lambda(x,r)\le(1-\lambda)r$;
\item $\gamma_\lambda(q_1,r)=(1-\lambda)r^2$ and $\gamma_\lambda(q_2,r)\asymp(1-\lambda)r$;
\item for every $x=\pi(ww'w'')$ with $w\in\mathbb{S}^{n_1}$ ending with $2$ or $3$, $w'\in\{0,1\}^{n_2}$ and $w''\in\mathbb{S}^\infty$ beginning by $2$ or $3$, for $r=2^{-(n_1+n_2)}$,
$$\gamma_\lambda\left(x,\lambda^{-1}r\right)=2^{n_1}(1-\lambda)\left(\lambda^{-1}r\right)^2\asymp 2^{-n_2}(1-\lambda)r\asymp 2^{-n_2}\gamma_\lambda\left(x,r\right).$$
Therefore, $\gamma_\lambda$ is never a weakly monotonic function, whatever $\lambda$ is -- and, of course, whatever $w_0/w_2$ is, since the capacity is independent of $\mu$.
\end{itemize}

In summary, $W$ is a scale (whose weak monotonicity follows the natural property of exit time function), while $\gamma_\lambda$ is only a quasi-scale.

\subsection{Heat kernel in the case $w_0=w_2$}

Recall that (cf.\ \cite{hk,L0}) for a p.c.f.\ fractal, $(\mathrm{E}_W)$ holds with scaling function $W(x,r)=r^\beta$ (where $\beta=\alpha+1\ge 2$ is called the walk dimension of that Dirichlet space, with $\alpha$ the Hausdorff dimension of that fractal), and the heat kernel exists, satisfying the sub-Gaussian estimate
$$\frac{1}{\mu(B(x,t^{1/\beta}))}1_{R(x,y)<t^{1/\beta}}\lesssim p_t(x,y)\lesssim\frac{1}{\mu(B(x,t^{1/\beta}))}\exp\left(-c\left(\frac{R(x,y)}{t^{1/\beta}}\right)^{\beta/(\beta-1)}\right)$$
for all $0<t\le 1$ and $x,y\in K$.

When we return to our present space $(K,R,\mu,\mathcal{E},\mathcal{F})$, such an estimate is impossible, even under weights $w_0=w_2$ (so that $\beta=3$). Actually, according to the framework of \cite{ghh+3,ghl15,L1}, the sub-Gaussian estimate implies $(\mathrm{EHI})$, contradicting \cite[Theorem 4.1]{HL}.

Below are some first attempts for a proper heat kernel estimate in the case $s_i=\frac{1}{2}$ and $w_i=\frac{1}{4}$ for $i\in\{0,1,2,3\}$.

\begin{lemma}\label{tail}
There exist constants $c,C>0$ such that for all $B=B(x_0,r)$ and $t>0$, on $\frac{1}{2}B$,
$$P_t1_{B^c}\le C\exp\left(-c\left(t^{-1}r^3\right)^{\frac{1}{2}}\right).$$
\end{lemma}

\begin{proof}
Recall by (\ref{imply}) that $(\mathrm{S}_W)$ holds with $W(x,r)=r^3$. The desired estimate follows from the same proof as of \cite[Theorem 3.4]{gh08}.
\end{proof}

Note that there is no counterpart under other weights $w_i$, since the iteration is valid only for scaling functions. See our previous work \cite[Lemma A.1]{CL} for detail how the scaling property (\ref{scal}) is used.

\begin{lemma}\label{nash}
The heat kernel $p_t(x,y)$ exists, and there is a constant $C>0$ such that
$$p_t(x,y)\le Ce^{t/4}t^{-4}$$
for all $x,y\in K$ and $t>0$.
\end{lemma}

\begin{proof}
For every $u\in\mathcal{F}$ and $n\in\mathbb{N}$, set $u_n(x)=\fint_{K_{w(x)}}u(y)d\mu(y)$, where $w(x)\in\mathbb{S}^n$ is taken such that the first address $\omega_x$ of $x$ begins with $w(x)$.

Note in particular that $x\in K_{w(x)}$. Hence for every $y\in K_{w(x)}$, we have
$$R(x,y)=2^{-n}R\left(F_{w(x)}^{-1}(x),F_{w(x)}^{-1}(y)\right)\le 2^{1-n},$$
yielding that
\begin{align*}
\|u-u_n\|_2^2=&\ \int_K(u-u_n)^2d\mu\le\sum\limits_{w\in\mathbb{S}^n}\int_{K_w}\left(u(x)-\fint_{K_w}u(y)d\mu(y)\right)^2d\mu(x)\\
=&\ \sum\limits_{w\in\mathbb{S}^n}\int_{K_w}\left(\fint_{K_w}(u(x)-u(y))d\mu(y)\right)^2d\mu(x)\\
\le&\ \sum\limits_{w\in\mathbb{S}^n}\int_{K_w}\left(\fint_{K_w}\sqrt{R(x,y)\mathcal{E}(u)}d\mu(y)\right)^2d\mu(x)\\
\le&\ \sum\limits_{w\in\mathbb{S}^n}\int_{K_w}2^{1-n}\mathcal{E}(u)d\mu=2^{1-n}\mathcal{E}(u).
\end{align*}

On the other hand, note that $\mu(K_w)=4^{-n}$ for every $w\in\mathbb{S}^n$. Therefore,
$$\|u_n\|_2^2=\int_K\left(\fint_{K_{w(x)}}u(y)d\mu(y)\right)^2d\mu(x)\le\int_K\left(4^n\int_Ku(y)d\mu(y)\right)^2d\mu(x)=2^{4n}\|u\|_1^2.$$

Combining these two inequalities, we obtain
$$\|u\|_2^2\le 2\|u-u_n\|_2^2+2\|u_n\|_2^2\le 2^{2-n}\mathcal{E}(u)+2^{1+4n}\|u\|_1^2.$$
Further, trivially $\|u\|_2^2\le s\|u\|_2^2$ for all $s\ge 1$. Thus by simple interpolation we see
$$\|u\|_2^2\le 2s\left(4\mathcal{E}(u)+\|u\|_2^2\right)+32s^{-4}\|u\|_1^2\quad\mbox{for all}\quad s>0.$$
Optimizing the right side over $s$, we see
$$\|u\|_2^2\le C\left(4\mathcal{E}(u)+\|u\|_2^2\right)^{\frac{4}{5}}\|u\|_1^{\frac{2}{5}}.$$
Comparing this inequality with \cite[Formula (2.2)]{cks}, the desired estimate follows.
\end{proof}

\begin{corollary}\label{UE}
For the case $w_0=w_2$, there exist constants $C,c>0$ such that for all $x,y\in K$ and $0<t\le 1$,
$$p_t(x,y)\le\frac{C}{t^4}\exp\left(-c\left(t^{-1}R(x,y)^3\right)^{\frac{1}{2}}\right).$$
\end{corollary}

This is definitely not a good estimate, since it contains no information on the measure. Improvement is possible, but an improvement to the on-diagonal upper estimate
\begin{equation}\tag*{$(\mathrm{DUE})$}
p_t(x,y)\le\frac{C}{\mu(B(x,t^{1/3}))}
\end{equation}
is impossible: since the strong Harnack inequality fails but a weak analogue holds, there must be some place where $p_t(x,\cdot)$ exceeds a given multiple of its on-diagonal average.

\begin{proof}
This follows from Lemmas \ref{tail} and \ref{nash}, using the well-known semigroup technique
\begin{align*}
p_t(x,y)=&\ \int_Kp_{\frac{t}{2}}(x,z)p_{\frac{t}{2}}(z,y)d\mu(z)\\
\le&\ \left(\int_{B(x,\frac{1}{2}R(x,y))^c}+\int_{B(y,\frac{1}{2}R(x,y))^c}\right)p_{\frac{t}{2}}(x,z)p_{\frac{t}{2}}(z,y)d\mu(z)\\
\le&\ \sup_{z\in K}p_{\frac{t}{2}}(z,y)\cdot P_{\frac{t}{2}}1_{B(x,\frac{1}{2}R(x,y))^c}(x)+\sup_{z\in K}p_{\frac{t}{2}}(x,z)\cdot P_{\frac{t}{2}}1_{B(y,\frac{1}{2}R(x,y))^c}(y)\\
\le&\ 2C\left(\frac{t}{2}\right)^{-4}\cdot C\exp\left(-c\left(\left(\frac{t}{2}\right)^{-1}\left(\frac{1}{2}R(x,y)\right)^3\right)^{\frac{1}{2}}\right)\\
=&\ \frac{32C^2}{t^4}\exp\left(-\frac{c}{2}\left(t^{-1}R(x,y)^3\right)^{\frac{1}{2}}\right),
\end{align*}
that is, the desired estimate with $32C^2$, $c/2$ replacing $C$, $c$ respectively.
\end{proof}

\begin{corollary}\label{LE}
For the case $w_0=w_2$, there exist constants $\varepsilon,\delta,\beta\in(0,1)$ such that for every ball $B=B(x_0,r)$ with $0<r\le 1$, for all $0<t\le\delta r^3$ and $x,y\in B(x_0,\delta t^{1/\beta})$,
$$p_t^B(x,y)\ge\frac{\varepsilon}{\mu(B(x_0,t^{1/3}))}.$$
\end{corollary}

This estimate is very weak. It can not be improved to the local lower estimate
\begin{equation}\tag*{$(\mathrm{LLE})$}
p_t^B(x,y)\ge\frac{C^{-1}1_{B(x,\varepsilon t^{1/3})}}{\mu(B(x,t^{1/3}))},
\end{equation}
since it yields $(\mathrm{EHI})$ by \cite[Corollary 2.5]{L1}. Meanwhile, it might be improved to $\beta>1$, so that by taking chains some off-diagonal lower estimate can be yielded.

\begin{proof}
This is modified from the proof of \cite[Proposition 5.3]{L1}, for which reason repeating details are omitted. Since $(\mathrm{S}_W)$ holds with $W(x,r)=r^3$, we see there exist $\varepsilon_1,\delta_1\in(0,1)$ such that
$$P_t^{B'}1_{B'}\ge\varepsilon_1\quad\mbox{on}\quad\delta_1B'$$
for all $0<t\le\delta_1r^3$ and $B'=B(x_0,(t/\delta_1)^{1/3})$. Say $\varepsilon_1,\delta_1\le\frac{1}{2}$. Combining the semigroup property, we have for all $x\in\delta_1B'$ that
\begin{equation}\label{diag-low}
p_t^{B'}(x,x)\ge\frac{1}{\mu(B')}\left(P_{\frac{t}{2}}^{B'}1_{B'}(x)\right)^2\ge\frac{\varepsilon_1^2}{\mu(B(x_0,t^{1/3}))}.
\end{equation}
On the other hand, by \cite[Theorem 4.2]{HL}, $(\mathrm{wEH}_1)$ holds. Thus by \cite[Lemma 6.2]{hy}, there exist constants $\gamma,\delta_2\in(0,1)$ such that
$$|u(x)-u(y)|\le C\left(\frac{d(x,y)}{r''}\right)^\gamma\|u\|_{L^\infty(B'')}\quad\mbox{for all}\quad x,y\in\delta_2B''$$
for every ball $B''=B(z,r'')$ and $u\in\mathcal{F}$ harmonic on $B''$. This implies (with the help of $(\mathrm{E}_W)$ and the Nash inequality) that, for all $r''<r'<(t/\delta_1)^{1/3}$,
$$\mathop{\mathrm{osc}}\limits_{\delta_2B(x,r'')}p_t^{B'}(x,\cdot)\le C\left(\frac{r''}{r'}\right)^\gamma\sup\limits_{B'\times B'}p_t^{B'}+C\frac{(r')^3}{t}\sup\limits_{y\in B'}p_{t/2}^{B'}(y,y)\le\frac{C}{t^4}\left\{\left(\frac{r''}{r'}\right)^\gamma+\frac{(r')^3}{t}\right\}$$
whenever $B(x,r'')\subset B'$. Optimizing the right side, we see
$$\mathop{\mathrm{osc}}\limits_{\delta_2B(x,r'')}p_t^{B'}(x,\cdot)\le\frac{C}{t^4}\left(\frac{(r'')^3}{t}\right)^{\frac{\gamma}{3+\gamma}}.$$
Combining (\ref{diag-low}), it follows that, for all $y\in B(x,r'')$,
$$p_t^{B'}(x,y)\ge p_t^{B'}(x,x)-\mathop{\mathrm{osc}}\limits_{\delta_2B(x,r'')}p_t^{B'}(x,\cdot)\ge\frac{\varepsilon_1^2}{\mu(B(x_0,t^{1/3}))}-\frac{C}{t^4}\left(\frac{(r'')^3}{t}\right)^{\frac{\gamma}{3+\gamma}}.$$
Recall that $C^{-1}r^2\le\mu(B(x_0,r))\le Cr$ for all $r\in(0,1]$. Therefore, by setting
$$r''=\left(\varepsilon_2t^{4-\frac{1}{3}}\right)^{\frac{3+\gamma}{3\gamma}}t^{1/3}=:\varepsilon'_2t^{1/\beta}$$
(where $\beta=\frac{9\gamma}{14\gamma+33}$), we see for all $x,y\in B(x_0,\delta_1\delta_2r'')$ that
$$p_t^B(x,y)\ge p_t^{B'}(x,y)\ge\frac{\varepsilon_1^2/2}{\mu(B(x_0,t^{1/3}))}.$$
Hence the desired result is proved with $\varepsilon=\varepsilon_1^2/2$ and $\delta=\delta_1\delta_2\varepsilon'_2$.
\end{proof}

\end{document}